# Newton polygons and formal groups: Conjectures by Manin and Grothendieck

By Frans Oort

## Introduction

We consider $p$-divisible groups (also called Barsotti-Tate groups) in characteristic $p$, their deformations, and we draw some conclusions.

For such a group we can define its Newton polygon (abbreviated NP). This is invariant under isogeny. For an abelian variety (in characteristic $p$) the Newton polygon of its $p$-divisible group is "symmetric".

In 1963 Manin conjectured that conversely any symmetric Newton polygon is "algebroid"; i.e., it is the Newton polygon of an abelian variety. This conjecture was shown to be true and was proved with the help of the "Honda-Serre-Tate theory". We give another proof in Section 5.

Grothendieck showed that Newton polygons "go up" under specialization: no point of the Newton polygon of a closed fiber in a family is below the Newton polygon of the generic fiber. In 1970 Grothendieck conjectured the converse: any pair of comparable Newton polygons appear for the generic and special fiber of a family. This was extended by Koblitz in 1975 to a conjecture about a sequence of comparable Newton polygons. In Section 6 we show these conjectures to be true.

These results are obtained by deforming the most special abelian varieties or $p$-divisible groups we can think of. In describing deformations we use the theory of *displays*; this was proposed by Mumford, and has been developed in [17], [18], and recently elaborated in [32] and [33]; also see [11], [31].

Having described a deformation we like to read off the Newton polygon of the generic fiber. In most cases it is difficult to determine the Newton polygon from the matrix defined by $F$ on a basis for the (deformed) Dieudonné module. In general I have no procedure to do this (e.g. in case we deform away from a formal group where the Dieudonné module is not generated by one element). However in the special case we consider here, $a(G_0) = 1$, a noncommutative version of the theorem of Cayley-Hamilton ("every matrix satisfies its own characteristic polynomial") suffices for our purposes. This enables us to find equations singling out specific Newton polygons (in that coordinate system NP-strata are linear subspaces). Once this is done (Sections 2 and 3), conjectures by Manin, Grothendieck and Koblitz follow easily (Sections 5 and 6).

Many people have patiently listened to me in discussions about this topic. Especially I mention, and I thank: Ching-Li Chai, Johan de Jong and Ben Moonen for sharing their time and interests with me.

## 1. Definitions and notation

Throughout the paper we fix a prime number $p$.

(1.1) For a *commutative* finite group scheme $N \to S$ we denote by $N^D$ its Cartier dual; see [22, I.2], [27, $\mathrm{VII}_A$.3], [30, 2.4]. It can be characterized functorially by:

$$N^D(T) = \mathrm{Hom}(N_T, \mathbb{G}_{m,T});$$

see [22, III.16], [27, $\mathrm{VII}_A$.3.3].

For an abelian scheme $\mathcal{X} \to S$ we denote its dual by $\mathcal{X}^t$. Note the duality theorem: for an isogeny $\varphi : \mathcal{X} \to \mathcal{Y}$ we canonically have an exact sequence

$$0 \to (\mathrm{Ker}(\varphi))^D \to \mathcal{Y}^t \xrightarrow{\varphi^t} \mathcal{X}^t \to 0;$$

see [22, III.19] and [16, III.15].

For $p$-divisible groups, see [28]. An abelian scheme $\mathcal{X} \to S$ of relative dimension $g$ defines a $p$-divisible group of height $2g$, which we denote by

$$\text{ind.lim.}\mathcal{X}[p^i] =: \mathcal{X}[p^\infty].$$

A *polarization* for an abelian scheme is an isogeny

$$\lambda : \mathcal{X} \to \mathcal{X}^t$$

which on each geometric fiber is defined by an ample divisor; see [15, 6.2]. Note that a polarization $\lambda$ is an isogeny which is symmetric in the sense that

$$(\lambda : \mathcal{X} \to \mathcal{X}^t)^t = \lambda$$

with the canonical identification $\mathcal{X} = \mathcal{X}^{tt}$; see [16, 21 Appl. III, p. 208]. A polarization is called a *principal* polarization if it is an isomorphism.

We say that $\lambda : G \to G^t$ is a quasi-polarization of a $p$-divisible group $G$ if it is a symmetric isogeny of $p$-divisible groups.

*From now on we work over a base scheme of characteristic $p$.*

(1.2) From [27, $\mathrm{VII}_A$.4], "Frobeniuseries", we recall: For every morphism $X \to S$, there is functorially a morphism

$$F : X \to X^{(p/S)}.$$

For a group scheme $G/S$ this is a homomorphism, and we write $G[F] := \mathrm{Ker}(F : G \to G^{(p/S)})$. For a *commutative* group scheme $G$ there is (functorially) a homomorphism

$$V : G^{(p/S)} \to G,$$

which moreover has the property that $V{\cdot}F = p = V{\cdot}F$.

(1.3) For a group scheme $G$ over a field $K$ of characteristic $p$ we define its $p$-rank $f(G) = f$ by: $\mathrm{Hom}(\mu_{p,k}, G_k) \cong (\mathbb{Z}/p)^f$; here $k$ is an algebraically closed field containing $K$, and $\mu_p := \mathbb{G}_m[p]$. Note that for an abelian variety $X$ this is the same as saying

$$X[p](k) \cong (\mathbb{Z}/p)^f.$$

We write $\alpha_p := \mathbb{G}_a[F]$; for a group scheme $G$ over a field $K$ of characteristic $p$,

$$a = a(G) := \dim_L(\mathrm{Hom}(\alpha_p, G \otimes_K L)),$$

where $L$ is a perfect field containing $K$ (this number does not depend on the choice of a perfect $L \supset K$).

(1.4) *Dieudonné modules.* In this paper we use the *covariant* theory. For a finite (*commutative*) group scheme $N$ of $p$-power rank over a (perfect) field $K$ there is a Dieudonné module $\mathbb{D}(N)$. This functor has the following properties:

- $N \mapsto \mathbb{D}(N)$ is an equivalence between the category of finite group schemes of $p$-power rank over $K$ and the category of Dieudonné modules of finite length as $W$-modules,
- if $\mathrm{rank}(N) = p^n$, then the length of $\mathbb{D}(N) = M$ equals $n$,
- $\mathbb{D}(F : N \to N^{(p)}) = (V : M \to M^{(p)})$,
- $\mathbb{D}(V : N^{(p)} \to N) = (F : M^{(p)} \to M)$.

In this ring $W[F, V]$ we have the relations $V{\cdot}F = p = F{\cdot}V$, and $F{\cdot}a = a^\sigma{\cdot}F$ and $V{\cdot}a^\sigma = a{\cdot}V$, for $a \in W = W_\infty(K)$; this ring is noncommutative if and only if $K \neq \mathbb{F}_p$.

For a $p$-divisible group $G$ of height $h$ over a perfect field $K$ there is a Dieudonné module $\mathbb{D}(G)$. This module is free of rank $h$ over $W$. If moreover the $p$-divisible group is a formal group, this module is over the ring $W[F][[V]]$. If $X$ is an abelian variety, we shall write $\mathbb{D}(X) = \mathbb{D}(X[p^\infty])$.

The dimension of $G$ is $d$ if $G[F]$ is of rank equal to $p^d$, and this is the case if and only if $\dim_K(\mathbb{D}(G)/V\mathbb{D}(G)) = d$. For a perfect field $K$ we have $a(G) = \dim_K(M/(FM + VM))$, where $M = \mathbb{D}(G)$.

(1.5) In [14], contravariant Dieudonné module theory was used. By $G_{m,n}$ we denote a $p$-divisible group defined over $\mathbb{F}_p$, given by: $G_{1,0} = \mathbb{G}_m[p^\infty]$ and $G_{0,1}$ for its dual, i.e. $G_{0,1}[p^i]$ is the constant group scheme $\mathbb{Z}/p^i$ over $\mathbb{F}_p$, and for coprime $m, n \in \mathbb{Z}_{>0}$ we define the (formal) $p$-divisible group $G_{m,n}$ by the covariant Dieudonné module

$$\mathbb{D}(G_{m,n}) = W[[F,V]]/W[[F,V]]{\cdot}(F^m - V^n)$$

(in the contravariant theory the exponents are interchanged; see [14, p. 35]). Note that the dimension of $G_{m,n}$ equals $m$, that

$$(G_{m,n})^t = G_{n,m},$$

hence the dual of $G_{m,n}$ has dimension $n$. The notation $G_{m,n}$ used in [14] and used here denote the same $p$-divisible group (but the Dieudonné module of it in [14] differs from the one used here).

By the Dieudonné-Manin classification, see [14] and [2], we find: if $G$ is a $p$-divisible group over a field $K$, there is a finite set of pairs $\{(m_i, n_i) \mid i \in I\}$ with $m_i \geq 0$ and $n_i \geq 0$ and $m_i$ and $n_i$ relatively prime, such that there is an isogeny

$$G \otimes k \quad \sim \quad \sum_{i \in I} \; G_{m_i,n_i},$$

where $k$ is an algebraically closed field containing $K$. This set of pairs is called the *formal isogeny type* of this $p$-divisible group $G$.

(1.6) *Notation* (*the Newton polygon*). The formal isogeny type of a $p$-divisible group is encoded in the concept of a Newton polygon, abbreviated NP. We write $\mathcal{N}(G)$ for the Newton polygon of $G$; each of the summands $G_{m,n}$ gives a slope $\lambda = n/(m+n)$ with multiplicity $(m+n)$; arranged in nondecreasing order this gives a polygon which has the following properties (the definition of a Newton polygon):

- The polygon starts at $(0,0)$ and ends at $(h, h-d)$ for a $p$-divisible group of height $h$ and dimension $d$,
- each slope $\lambda \in \mathbb{Q}$ has the property $0 \leq \lambda \leq 1$,
- the polygon is lower convex, and
- its break-points have integral coordinates.

By "lower convex" we mean that it is the graph of a convex piecewise-linear function on the interval $[0, h]$. A Newton polygon determines, and is determined by, its set of slopes $0 \leq \lambda_1 \leq \cdots \leq \lambda_h \leq 1$. It is called *symmetric* if $\lambda_i = 1 - \lambda_{h-i+1}$ for $1 \leq i \leq h$.

(1.7) *The formal isogeny type of an abelian variety $X$ is* symmetric; *i.e., it can be written*

$$X[p^\infty] \otimes k \quad \sim \quad f{\cdot}(G_{1,0} \oplus G_{0,1}) \quad \oplus \quad s{\cdot}G_{1,1} \quad \oplus \sum_j \ (G_{m_j,n_j} \oplus G_{n_j,m_j}).$$

Indeed, an abelian variety has a polarization; by the duality theorem this implies that $X[p^\infty]$ is isogenous with its dual. □

The converse of this statement is called the "Manin conjecture"; see (5.1).

(1.8) We say that a Newton polygon $\beta$ is *lying above* $\gamma$, notation

$$\beta \quad \prec \quad \gamma,$$

if $\beta$ and $\gamma$ have the same end points, and if no point of $\beta$ is strictly below $\gamma$ (!! note the reverse order). The Newton polygon consisting only of slopes $\frac{1}{2}$ is called the *supersingular* one, notation $\sigma$; this is a symmetric Newton polygon. The Newton polygon belonging to $d{\cdot}G_{1,0} \oplus c{\cdot}G_{0,1})$ is called the *ordinary* one, denoted by $\rho = \rho_{d,c}$. Note that any symmetric $\xi$ satisfies $\sigma \prec \xi \prec \rho = \rho_{g,g}$. More generally, a Newton polygon of height $h$ and dimension $d$ is between the straight line, the Newton polygon of $G_{d,d-h}$ and the Newton polygon of $dG_{1,0} + (h-d)G_{0,1}$.

(1.9) *Displays over a field. In this section we work over a perfect field* $K \supset \mathbb{F}_p$. Covariant Dieudonné module theory over a perfect field is a special case of the theory of displays. Consider a $p$-divisible group $G$ of height $h$ over $K$, and its Dieudonné-module $\mathbb{D}(G) = M$. We choose a $W$-base

$$\{e_1 = X_1, \ldots, e_d = X_d, \quad e_{d+1} = Y_1, \ldots, e_h = Y_c\}$$

for $M$ such that $Y_1, \ldots, Y_c \in VM$; on this base the structure of the Dieudonné module is written as:

$$\begin{aligned} Fe_j &= \sum_{i=1}^{h} a_{ij} e_i, && 1 \leq j \leq d, \\ e_j &= V\left(\sum_{i=1}^{h} a_{ij} e_i\right), && d < j \leq h. \end{aligned}$$

We have written the module in displayed form; see [17, §1], [18, §0], and [32]. Now

$$(a_{i,j} \mid 1 \leq i, j \leq h) = \begin{pmatrix} A & B \\ C & D \end{pmatrix}.$$

This matrix, denoted by $(a)$, will be called the matrix of the display. Note that in this case the $\sigma$-linear map $F$ is given on this base by the matrix

$$\begin{pmatrix} A & pB \\ C & pD \end{pmatrix},$$

where

$$\begin{aligned} A &= (a_{ij} \mid 1 \leq i,j \leq d), \quad B = (a_{ij} \mid 1 \leq i \leq d < j \leq h), \\ C &= (a_{ij} \mid 1 \leq j \leq d < i \leq h,), \quad D = (a_{ij} \mid d < i,j \leq h). \end{aligned}$$

Where the display-matrix is symbolically denoted by $(a)$, we write $(pa)$ symbolically for the associated $F$-matrix (it is clear what is meant as soon as $d$ is given). Note that the induced maps

$$F : M/VM \longrightarrow FM/pM, \quad \frac{F}{p} : VM/pM \longrightarrow M/FM$$

are bijective, hence the matrix $(a_{i,j} \mid 1 \leq i,j \leq h)$ has as determinant a unit in $W$; let its inverse be $(b_{i,j})$, written in block form as

$$(b_{i,j}) = \begin{pmatrix} E & G \\ H & J \end{pmatrix}$$

(in [32], the block matrix called $J$ here is denoted by $B$). Working over a perfect field $K$, with $W := W_\infty(K)$, we write

$$\tau := \sigma^{-1} : \quad W \longrightarrow W.$$

Then the map $V : M \to M$ on the given basis has as matrix

$$\begin{pmatrix} pE & pG \\ H & J \end{pmatrix}^{\tau},$$

called the $V$-matrix.

Note that the $p$-divisible group is a formal group if and only if the operation $V$ on its covariant Dieudonné module is topologically nilpotent (note that $\mathbb{D}(F) = V$). We remark that the Dieudonné module of $G[p]/G[F]$ corresponds with $\oplus_{i=d+1}^{h} K{\cdot}e_i$; hence we see that $G$ is a formal $p$-divisible group if and only if the matrix $J$ mod $p$ is nilpotent in the $\tau$-linear sense (also see [32, p. 6]).

We write $Q := VM$, $\; T := \oplus_{i=1}^{d} W{\cdot}e_i$, $\; L := \oplus_{i=d+1}^{h} W{\cdot}e_i$, and note that

$$F \oplus V^{-1} : T \oplus L \to M,$$

as given above by the transformation formulas on a basis, is a $\sigma$-linear bijective map (we follow [32] for this notation).

(1.10) Suppose $(G, \lambda)$ is a $p$-divisible group with a principal quasi-polarization over a perfect field $K$. The quasi-polarization can be given on the Dieudonné module $M = \mathbb{D}(G)$ by a skew perfect pairing

$$\langle \ , \ \rangle : M \quad \times \quad M \longrightarrow W$$

which satisfies

$$\langle Fa, b\rangle = \langle a, Vb\rangle^{\sigma}, \qquad \text{for all } a,\ b \in M;$$

see [20, p. 83]. Then we can choose a *symplectic* base $\{X_1, \ldots, X_d, Y_1, \ldots, Y_d\}$ for the Dieudonné module $\mathbb{D}(G) = M$; i.e., the polarization is given by a skew bilinear form which on this base is given by:

$$\langle X_i, Y_j\rangle = \delta_{ij}, \quad \langle X_i, X_j\rangle = 0 = \langle Y_i, Y_j\rangle.$$

Note that if the module is in displayed form on a symplectic base as above, with the display-matrices $(a)$ and $(b)$ as given above, then not only do we have $(a)\cdot(b) = \mathbf{1} = (b)\cdot(a)$ but also

$$(a_{i,j}) = \begin{pmatrix} A & B \\ C & D \end{pmatrix}, \quad (b_{i,j}) = \begin{pmatrix} E & G \\ H & J \end{pmatrix} = \begin{pmatrix} D^t & -B^t \\ -C^t & A^t \end{pmatrix},$$

where $A^t$ is the transpose of the matrix $A$.

(1.11) *Displays.* In order to describe deformations we choose a complete Noetherian local ring $R$ with perfect residue class field $K$ (a complete Noetherian local ring is excellent). We assume $p{\cdot}1 = 0 \in R$. In this case the theory of displays as described in [32] gives *an equivalence of categories between the category of displays over $R$ and the category of formal $p$-divisible groups over $R$.* We refer to [18], [32] and [33] for further details and will describe deformations by constructing a display.

(1.12) *Deformations of formal $p$-divisible groups.* Displays can be applied as follows. Consider a *formal* $p$-divisible group $G_0$ over a perfect field $K$, and suppose we have written its Dieudonné module and base $\{e_1, \cdots, e_h\} = \{X_1, \ldots, X_d, Y_1, \ldots, Y_c\}$ in displayed form as above. We write $h$ for the height of $G_0$, and $d$ for its dimension, $c := h - d$. Choose $R$ as above, let $t_{r,s}$ be elements in the maximal ideal of $R$, with $1 \le r \le d < s \le h$, and let

$$T_{r,s} = (t_{r,s}, 0, \cdots) \in W(R)$$

be their Teichmüller lifts. We define a display over $R$ by considering a base

$$\{e_1 = X_1, \ldots, e_d = X_d, \quad e_{d+1} = Y_1, \ldots, e_h = Y_c\},$$

and

$$FX_j = \sum_{i=1}^{d} a_{ij} X_i + \sum_{v=1}^{h-d} a_{d+v,j}\Big(Y_v + \sum_{r=1}^{d} T_{r,v+d} X_r\Big), \quad 1 \le j \le d,$$

$$Y_t = V\left(\sum_{i=1}^{h} a_{i,d+t} X_i + \sum_{v=1}^{h-d} a_{d+v,d+t}\Big(Y_v + \sum_{r=1}^{d} T_{r,v+d} X_r\Big)\right), \quad 1 \le t \le c = h{-}d.$$

Note that this corresponds with the matrix

$$\begin{pmatrix} A+TC & pB+pTD \\ C & pD \end{pmatrix}, \qquad T = \begin{pmatrix} T_{1,d+1} & \cdots & T_{1,h} \\ \vdots & & \vdots \\ T_{d,d+1} & \cdots & T_{d,h} \end{pmatrix},$$

which gives the map $F$. Now if for display-matrices,

$$\begin{pmatrix} A & B \\ C & D \end{pmatrix}^{-1} = \begin{pmatrix} E & G \\ H & J \end{pmatrix},$$

then

$$\begin{pmatrix} A+TC & B+TD \\ C & D \end{pmatrix}^{-1} = \begin{pmatrix} E & G-ET \\ H & J-HT \end{pmatrix}.$$

(1.13) *Let $K$ be a perfect field, let $R$ be a complete Noetherian local ring with residue class field $K$, and let $G_0$ be a formal $p$-divisible group over $K$. The formulas above define a display over $W$, the ring of Witt vectors over $R$* (in the sense of [32]). *Hence these formulas define a deformation $G \to \mathrm{Spec}(R)$ of $G_0$.*

In fact, we have written the module in displayed form. We write $\epsilon : W \to K$ for the residue class map. The matrix $\epsilon(J - HT)$ is $\sigma$-linear nilpotent (still, after deforming); hence this gives the condition necessary for a "display" in the sense of [32]. By the theory of Dieudonné modules we conclude that this defines a formal $p$-divisible group $G \to \mathrm{Spec}(R)$; see [17], [18], [32], [33].

(1.14) *Remark.* One can show that the deformation just given is the universal deformation of $G_0$ in equal characteristic $p$, by taking the elements $t_{r,s}$ as parameters,

$$R := K[[t_{r,s} \mid 1 \le r \le d < s \le h]].$$

(1.15) *Deformations of principally quasi-polarized formal $p$-divisible groups.* Suppose moreover the $p$-divisible group $G_0$ has a principal quasi-polarization $\lambda_0$ and let the base $\{X_1, \ldots, Y_d\}$ be in symplectic form (in this case $c = h-d = d$). Assume moreover that

$$t_{r,s} = t_{s-d,r+d} \in R,$$

and the displayed form above defines a deformation $(G, \lambda)$ as a quasi-polarized formal $p$-divisible group of $(G_0, \lambda_0)$ (see [18], [32], [33]). After renumbering: $x_{i,j} = t_{i,j+d}$, with $d = g = c = h/2$, we have the familiar equations $x_{i,j} = x_{j,i}$, $1 \le i, j \le g$.

## 2. Cayley-Hamilton

We denote by $k = \overline{k} \supset \mathbb{F}_p$ an algebraically closed field.

(2.1) *Definition.* We consider matrices which can appear as $F$-matrices associated with a display. Let $d, c \in \mathbb{Z}_{\geq 0}$, and $h = d + c$. Let $W$ be a ring. We say that a display-matrix $(a_{i,j})$ of size $h \times h$ is in *normal form* form over $W$ if the $F$-matrix is of the following form:

$$(\mathcal{F}) \qquad \begin{pmatrix} 0 & 0 & \cdots & 0 & a_{1d} & pa_{1,d+1} & \cdots & \cdots & \cdots & pa_{1,h} \\ 1 & 0 & \cdots & 0 & a_{2d} & \cdots & & pa_{i,j} & & \cdots \\ 0 & 1 & \cdots & 0 & a_{3d} & & & 1 \leq i \leq d & & \\ \vdots & \vdots & \ddots & \ddots & \vdots & & & d \leq j \leq h & & \\ 0 & 0 & \cdots & 1 & a_{dd} & pa_{d,d+1} & \cdots & \cdots & \cdots & pa_{d,h} \\ & & & & & & & & & \\ 0 & \cdots & \cdots & 0 & 1 & 0 & \cdots & \cdots & \cdots & 0 \\ 0 & \cdots & & \cdots & 0 & p & 0 & \cdots & \cdots & 0 \\ 0 & \cdots & & \cdots & 0 & 0 & p & 0 & \cdots & 0 \\ 0 & \cdots & & \cdots & 0 & 0 & 0 & \ddots & 0 & 0 \\ 0 & \cdots & & \cdots & 0 & 0 & \cdots & \cdots & p & 0 \end{pmatrix},$$

$a_{i,j} \in W, \quad a_{1,h} \in W^*$; i.e. it consists of blocks of sizes ($d$ or $c$) $\times$ ($d$ or $c$); in the left hand upper corner, which is of size $d \times d$, there are entries in the last column, named $a_{i,d}$, and the entries immediately below the diagonal are equal to 1; the left and lower block has only one element unequal to zero, and it is 1; the right hand upper corner is unspecified, entries are called $pa_{i,j}$; the right hand lower corner, which is of size $c \times c$, has only entries immediately below the diagonal, and they are all equal to $p$.

Note that if a Dieudonné module is defined by a matrix in displayed normal form then either the $p$-rank $f$ is maximal, $f = d$, and this happens if and only if $a_{1,d}$ is not divisible by $p$, or $f < d$, and in that case $a = 1$. The $p$-rank is zero if and only if $a_{i,d} \equiv 0 \pmod{p}$, for all $1 \leq i \leq d$.

(2.2) Lemma. *Let $M$ be the Dieudonné module of a $p$-divisible group $G$ over $k$ with $f(G) = 0$. Suppose $a(G) = 1$. Then there exists a $W$-basis for $M$ on which $F$ has a matrix which is in normal form.*

In this case the entries $a_{1,d}, \ldots, a_{d,d}$ are divisible by $p$ and can be chosen to equal zero.

The proof is easy and is omitted, but we do give the proof of the following lemma which is slightly more involved.

(2.3) Lemma. *Let $k$ be an algebraically closed field, and let $M$ be the Dieudonné module of a local, principally quasi-polarized formal $p$-divisible group $G$ over $k$ with $f(G) = 0$ (hence $G$ is of "local-local type"). Write $d$ for its dimension, and $h = 2d$ for its height. Suppose $a(M) = 1$. Then there exists a* symplectic *$W$-basis for $M$ on which the matrix of $F$ is in normal form; see* (2.1).

*Proof.* For every $a \in \mathbb{Z}_{>0}$ we shall choose an element $X^{(a)}$, with $X^{(a)} \notin FM + VM$, with $X^{(a+1)} - X^{(a)} \in p^a M$.

*First step.* We choose $X \in M$ with $X \notin FM + VM$. Let $\langle X, F^d X\rangle = \beta \in W$, and write $\beta \bmod pW = \beta_0 \in k$. We note that $\beta_0 \neq 0$; in fact, write $FM_1 \cap VM_1 = F^d M_1 = V^d M_1$; this is of dimension one over $k$; note that the pairing induces a perfect pairing between $M/(FM + VM)$ and $F^d M_1$. Choose $\lambda \in W$, with $\lambda_0 \in k = W/pW$ such that $\lambda_0^{p^d+1}{\cdot}\beta_0 = 1$. Note that $\langle \lambda X, F^d \lambda X\rangle = \lambda \cdot \lambda^{p^d} \cdot \beta$; hence $\langle \lambda X, F^d \lambda X\rangle \equiv 1 \bmod pW$.

We change notation, choosing a new $X$ instead of the old $\lambda \cdot X$, and concluding that for the new $X$,

$$\langle X, F^d X\rangle \equiv 1 \bmod pW.$$

Suppose

$$\langle X, F^{s-1} X\rangle = b_s, \quad 2 \leq s \leq d.$$

Choosing

$$X' = X + \sum_{s=2}^{d} b_s{\cdot}Y_s,$$

we see that

$$\langle X', F^{s-1} X'\rangle \equiv 0 \pmod{pW},$$

and still $\langle X', F^d X'\rangle \equiv 1$.

We call this new $X'$ finally $X^{(1)}$, defining $X_s = F^{s-1} X^{(1)}$ for $1 \leq s < d$, and $Y_j = -V^{d-j+1} X^{(1)}$ for $1 \leq j \leq d$. The set $\{X^{(1)} = X_1, \ldots, Y_d\}$ is a $W$-basis for $M$, which is symplectic modulo $pW$. On this basis the matrix of $F$ is congruent to $(\mathcal{F})$ mod $pW$.

We introduce some notation to be used in this proof. We say that $B = \{X_1, \ldots, X_d, Y_1, \ldots, Y_d\}$ is an n-basis, if $X = X^{(a)} \in W$, with $X \notin FM + VM$, and $X_s := F^{s-1} X$ for $1 < s \leq d$, and $Y_j \equiv -V^{d-j+1} X \bmod pW$ for $1 \leq j \leq d$; note that an n-basis indeed is a $W$-basis for $M$. We say that an n-basis is (sa) if it is simplectic modulo $p^a W$. We say that an n-basis is (Fa) if the matrix of $F$ on this basis is in normal form; see (2.1), modulo $p^a W$.

For $a \in \mathbb{Z}_{>0}$ we study the statement:

$(I_a)$: *For* $X = X^{(a)}$, *assume that there exists an* n-*basis* $B = B^{(a)}$ *which is* (sa) *and* (Fa).

The construction in the first step makes a choice such that $(I_1)$ is satisfied. This ends the first step.

We assume for $a \geq 1$ that $(I_a)$ is satisfied, and construct a basis which satisfies $(I_{a+1})$. In order to formulate the proof of the induction step we assume that $a \in \mathbb{Z}_{>0}$ is fixed, and write $\equiv$ for equivalence modulo $p^aW$, and $\approx$ for equivalence modulo $p^{a+1}W$. In order to prove the induction step $(I_a) \mapsto (I_{a+1})$ we study the condition:

$(J_k)$ Here $a \in \mathbb{Z}_{>0}$ is fixed, and $1 \leq k \leq d$. *Assume there is an* n-*basis satisfying* $(I_a)$, *and satisfying the property that* $\{X_1, \ldots, X_d, Y_1, \ldots, Y_k\}$ *is* $(ps(a+1))$. This last property means that pairs of elements appearing in this set satisfy the symplectic properties modulo $p^{a+1}W$, *and* $Y_1, \ldots, Y_{k-1}$ *transform as prescribed by some* $(\mathcal{F})$ mod $p^{a+1}W$.

*Second Step. We choose a basis satisfying* $(J_1)$. By induction we know that $\langle X, F^dX\rangle \equiv 1 + p^a{\cdot}\beta$, with $\beta \in W$; write $\beta_0 = \beta \mod pW$. Choose $\lambda \in W$ such that $\lambda_0^{p^d} + \lambda_0 + \beta_0 = 0$, and replace $X$ by $(1 + p^a{\cdot}\lambda)X$. For this new $X$ we have achieved $\langle X, F^dX\rangle \approx 1$. Choose $b_s \in W$ such that for this new $X$ we have: $\langle X, F^{s-1}X\rangle \approx p^a{\cdot}b_s$, for $2 \leq s \leq d$. Now, choose

$$X' := X + \sum_{s=2}^{d} b_s{\cdot}p^a{\cdot}Y_s,$$

and then

$$\langle X', F^{s-1}X'\rangle \equiv 0 \pmod{p^{a+1}W}.$$

We call this $X'$ finally $X$. The new $X, FX, \ldots, Y_1 = F^dX$ and the old $Y_2, \ldots, Y_d$ satisfy the condition $(J_1)$. This ends the second step.

*Third Step.* Here $a \in \mathbb{Z}_{>0}$ is fixed. With $1 \leq k < d$, we assume $(J_k)$ satisfied, and construct and prove: *There exists a basis satisfying* $(J_k)$ *such that* $\langle Y_k, VX\rangle \equiv 0 \pmod{p^{a+2}W}$.

Next we construct and prove $(J_{k+1})$.

There exists an element $\xi \in p^{a+2}M$. Choose a new $X$ instead of the old $X+\xi$, construct $X_1, \ldots, X_d, Y_1, \ldots, Y_k$ which together with the old $Y_{k+1}, \ldots, Y_d$ satisfy $(J_k)$ and the properties: $\langle X_i, X_j\rangle \equiv 0 \mod p^{a+2}$, for all $1 \leq i, j \leq d$, and $\langle FY_k, X_i\rangle \equiv 0 \mod p^{a+2}$ for all $i$ and $\langle FY_k, X_{k+1}\rangle \equiv 1 \mod p^{a+2}$. We choose this new n-basis and write

$$FY_k = pY' + px \quad \text{with} \quad Y' \in W{\cdot}Y_1 + \cdots W{\cdot}Y_d, \quad \text{and} \quad x \in W{\cdot}X_1 + \cdots W{\cdot}X_d.$$

Let $\langle FY_k, Y_j\rangle = c_j$ for $1 \leq j \leq k$; note that $c_j \in p^{a+1}W$. Define the new

$$Y_{k+1} := Y' - \sum_{1\leq j\leq k} c_j X_j; \quad \text{then} \quad FY_k = p{\cdot}Y_{k+1} + p\left(x + \sum_{1\leq j\leq k} c_j X_j\right).$$

The new $X_1, \ldots, X_d, Y_1, \ldots, Y_{k+1}$ together with the old $Y_{k+2}, \ldots, Y_d$ satisfy the condition $(J_{k+1})$. This ends the third step.

For fixed $a \in \mathbb{Z}_{>0}$, we start induction $(J_1)$ by the second step, and then by induction in the third step $(J_k) \mapsto (J_{k+1})$ we show that $(J_d)$ is satisfied. Note that $(I_a) + (J_d) = (I_{a+1})$; thus we have proved that $(I_a) \Rightarrow (I_{a+1})$. By the first step we can start induction: $I_1$ is satisfied. Hence induction shows that all steps $(I_a)$ for $a \in \mathbb{Z}_{>0}$ are satisfied. Moreover the bases $B^{(a)} = \{X_1^{(a)}, \ldots X_d^{(a)}, Y_1^{(a)}, \ldots, Y_d^{(a)}\}$ constructed satisfy $B^{(a+1)} \equiv B^{(a)}$ (mod $p^a M$). Hence this process converges to a $W$-basis for $M$; by construction, on this basis the matrix of $F$ is in normal form. This ends the proof of the lemma. □

(2.4) *How to determine* $\mathcal{N}(G)$? Consider the Dieudonné module of a $p$-divisible group (or of an abelian variety), by writing down the matrix of $F$ relative to some $W$-basis. In general it is difficult to see directly from that matrix what the Newton polygon of the $p$-divisible group studied is.

There are general results which enable us to compute the Newton polygon from a given displayed form. Here is an example: Nygaard proved in [19, Th. 1.2, p. 84], a general result, which e.g. for $g = 3$ gives the following: Suppose there is an abelian variety of dimension 3, and let $F$ be the action of Frobenius on its Dieudonné module; this abelian variety is supersingular if and only if $p^3|F^8$. For explicit computations this does not look attractive. Also see [9]: "slope estimates".

(2.5) *Remark.* We can compute the $p$-adic values of the eigenvalues of the matrix. Note however that if we change the basis, this $\sigma$-linear map gives a matrix on the new basis in the $\sigma$-linear way. In [9, pp. 123/124], we find an example by B. H. Gross of a $2 \times 2$ matrix which, on one basis has $p$-adic values of the eigenvalues equal to $\frac{1}{2}$, while after a change of basis these $p$-adic values equal 0 and 1. We see that the change of basis can change the $p$-adic values of the eigenvalues of this matrix. Thus we have the question: how can we determine the Newton polygon from the matrix (say of $F$ on the Dieudonné module)?

In this section we show how in case the matrix is in *normal form* the Newton polygon can be read off easily:

(2.6) LEMMA (of Cayley-Hamilton type). *Let $L$ be a field of characteristic $p$, let $W = W_\infty(L)$ be its ring of infinite Witt vectors. Let $G$ be a $p$-divisible group, with* $\dim(G) = d$, *and height$(G) = h$, with Dieudonné module $M$. Suppose there is a $W$-basis of $M$, such that the display-matrix $(a_{i,j})$ on this base gives an $F$-matrix in normal form as in* (2.1). *Now, $X = X_1 = e_1$ for the first base vector. Then for the expression*

$$P := \sum_{i=1}^{d} \sum_{j=d}^{h} p^{j-d} a_{i,j}^{\sigma^{h-j}} F^{h+i-j-1}, \quad F^h \cdot X = P \cdot X.$$

Note that we take powers of $F$ in the $\sigma$-linear sense, i.e. if the display matrix is $(a)$, i.e. $F$ is given by the matrix $(pa)$ as above,

$$F^n \quad \text{is given by the matrix} \quad (pa)\cdot(pa^\sigma)\cdot \dots \cdot(pa^{\sigma^{n-1}}).$$

The exponent $h+i-j-1$ runs from $0 = h+1-h-1$ to $h-1 = h+d-d-1$. Note that we do not claim that $P$ and $F^h$ have the same effect on all elements of $M$.

*Proof.* Note that $F^{i-1}e_1 = e_i$ for $i \leq d$.

CLAIM. *For $d \leq s < h$,*

$$F^s X = \left( \sum_{i=1}^{d} \sum_{j=d}^{s} F^{s-j} p^{j-d} a_{i,j} F^{i-1} \right) X + p^{s-d} e_{s+1}.$$

This is correct for $s = d$. The induction step from $s$ to $s+1 < h$ follows from $Fe_{s+1} = \left(\sum_{i=1}^{d} p\, a_{i,s+1} F^{i-1}\right) X + pe_{s+2}$. This proves the claim. Computing $F(F^{h-1}X)$ gives the desired formula. □

(2.7) PROPOSITION. *Let $k$ be an algebraically closed field of characteristic $p$, let $W = W_\infty(K)$ be its ring of infinite Witt vectors. Suppose $G$ is a $p$-divisible group over $k$ such that for its Dieudonné module the map $F$ is given by a matrix in normal form. Let $P$ be the polynomial given in the previous proposition. The Newton polygon $\mathcal{N}(G)$ of this $p$-divisible group equals the Newton polygon given by the polynomial $P$.*

*Proof.* Consider the $W[F]$-sub-module $M' \subset M$ generated by $X = e_1$. Note that $M'$ contains $X = e_1, e_2, \dots, e_d$. Also it contains $Fe_d$, which equals $e_{d+1}$ plus a linear combination of the previous ones; hence $e_{d+1} \in M'$. In the same way we see: $pe_{d+2} \in M'$, and $p^2 e_{d+3} \in M'$ and so on. This shows that $M' \subset M = \oplus_{i \leq h}$ and $W\cdot e_i$ is of finite index and that $M' = W[F]/W[F] \cdot(F^h - P)$. From this we see by the classification of $p$-divisible groups up to isogeny, that the result follows by [14, II.1]. Also see [2, pp. 82-84]. By [2, p. 82, Lemma 2] we conclude that the Newton polygon of $M'$ in case of the monic polynomial $F^h - \sum_0^m \ b_i F^{m-i}$ is given by the lower convex hull of the pairs $\{(i, v(b_i)) \mid i\}$. Hence the proposition is proved. □

(2.8) Corollary. *With notation as above, suppose that every element* $a_{i,j},\ \ 1 \le i \le d,\ \ d \le j \le h$, *is either equal to zero, or is a unit in* $W(k)$. *Let* $S$ *be the set of pairs* $(i,j)$ *with* $0 \le i \le d$ *and* $d \le j \le h$ *for which the corresponding element is nonzero*:

$$(i,j) \in S \quad \Longleftrightarrow \quad a_{i,j} \neq 0.$$

*Consider the image* $T$ *under*

$$S \to T \subset \mathbb{Z} \times \mathbb{Z} \quad \textit{given by} \quad (i,j) \mapsto (j+1-i, j-d).$$

*Then* $\mathcal{N}(G)$ *is the lower convex hull of the set* $T \subset \mathbb{Z} \times \mathbb{Z}$ *and the point* $(0,0)$; *note that* $a_{1,h} \in W^*$, *hence* $(h, h-d=c) \in T$.

This can be visualized in the following diagram (illustrating the case $d \le h-d$):

$$\begin{array}{ccccccccccc}
 & & & & & & & a_{d,h} & \cdots & & a_{1,h} \\
 & & & & & & \cdot & \cdots & & \cdot & \\
 & & & & & a_{d,2d+2} & \cdot & \cdot & \cdot & & \\
 & & & & a_{d,2d+1} & \cdot & \cdot & \cdot & a_{1,2d+1} & & \\
 & & & \cdot & \vdots & \vdots & \vdots & \cdot & & & \\
 & a_{d,d+1} & \cdots & a_{i,d+1} & \cdots & a_{2,d+1} & a_{1,d+1} & & & & \\
a_{d,d} & \cdots & a_{i,d} & \cdots & \cdots & a_{1,d} & & & & &
\end{array}$$

Here the element $a_{d,d}$ is in the plane with coordinates $(x=1, y=0)$ and $a_{1,h}$ has coordinates $(x=h, y=h-d=c)$. One erases the spots where $a_{i,j}=0$, and leaves the places where $a_{i,j}$ is a unit. The lower convex hull of these points and $(0,0)$ (and $(h, h-d)$) equals $\mathcal{N}(G)$.

*Proof.* This we achieve by writing out the Newton polygon of the polynomial $P$ in the Cayley-Hamilton lemma. □

## 3. Newton polygon strata

In this section $k = \overline{k} \supset \mathbb{F}_p$ is an algebraically closed field.

(3.1) We fix integers $h \ge d \ge 0$, and write $c := h-d$ considering Newton polygons ending at $(h,c)$. For such a Newton polygon $\beta$,

$$\Diamond(\beta) = \{(x,y) \in \mathbb{Z} \times \mathbb{Z} \mid y < c, \ \ y < x, \ \ (x,y) \prec \beta\};$$

here we denote by $(x, y) \prec \beta$ the property "$(x, y)$ is on or above $\beta$"; now,

$$\dim(\beta) = \#(\Diamond(\beta)).$$

Note that for $\rho = d{\cdot}(1,0) + c{\cdot}(0,1)$, $\dim(\rho) = dc$.

Suppose there is a formal $p$-divisible group $G_0$ over $k$ with $\mathcal{N}(G_0) = \gamma$ ending at $(h, c)$. We write $\mathcal{D} = \mathrm{Def}(G_0)$ for the universal deformation space in equal characteristic $p$. By this we mean the following. Formal deformation theory of $G_0$ is prorepresentable; we obtain a formal scheme $\mathrm{Spf}(A)$ and a prorepresenting family $\mathcal{G}' \to \mathrm{Spf}(A)$. By [5, Lemma 2.4.4, p. 23], we know that there is an equivalence of categories of $p$-divisible groups over $\mathrm{Spf}(A)$ respectively over Spec(A). We will say that $\mathcal{G} \to \mathrm{Spec}(A) = \mathcal{D} = \mathrm{Def}(G_0)$ is the universal deformation of $G_0$ if the corresponding $\mathcal{G}' \to \mathrm{Spf}(A) = \mathcal{D}^\wedge$ prorepresents the deformation functor.

A theorem by Grothendieck and Katz, see [9, Th. 2.3.1, p. 143], says that for any family $\mathcal{G} \to S$ of $p$-divisible groups over a scheme $S$ in characteristic $p$, and for any Newton polygon $\beta$ there is a unique maximal closed, reduced subscheme $W \subset S$ containing all points $s$ at which the fiber has a Newton polygon equal to or lying above $\beta$:

$$s \in W \quad \Longleftrightarrow \quad \mathcal{N}(\mathcal{G}_s) \prec \beta.$$

This set will be denoted by

$$\mathcal{W}_\beta(\mathcal{G} \to S) \subset S.$$

For every Newton polygon $\beta$ with $\beta \succ \gamma = \mathcal{N}(G_0)$ we define $V_\beta \subset \mathcal{D} = \mathrm{Def}(G_0)$ as the maximal closed, reduced subscheme carrying all fibers with Newton polygon equal to or above $\beta$ in the universal deformation space of $G_0$. Note that $V_\rho = \mathcal{D}$, with $\rho = d(1,0) + c(0,1)$.

In case of a family $(\mathcal{G}, \lambda) \to S = \mathrm{Spec}(R)$ of quasi-polarized $p$-divisible groups there is an analogous notion, and for a symmetric Newton polygon $\xi$ we write

$$\mathcal{W}_\xi(\mathcal{G} \to \mathrm{Def}(G_0, \lambda_0)) =: W_\xi \ \subset \ \mathrm{Def}(G_0, \lambda_0).$$

(3.2) Theorem (Newton polygon-strata for formal groups). *Suppose* $a(G_0) \le 1$. *For every* $\beta \succ \gamma = \mathcal{N}(G_0)$, $\dim(V_\beta) = \dim(\beta)$. *The strata* $V_\beta$ *are nested as given by the partial ordering on Newton polygons*; *i.e.*,

$$V_\beta \subset V_\delta \Longleftrightarrow \Diamond(\beta) \subset \Diamond(\delta) \Longleftrightarrow \beta \prec \delta.$$

*Generically on* $V_\beta$ *the fibers have Newton polygon equal to* $\beta$. *There is a coordinate system on* $\mathcal{D}$ *in which all* $V_\beta$ *are linear subspaces.*

*Proof.* At first we choose a coordinate system for $\mathcal{D} = \mathrm{Def}(G_0)$ where

$$G_0 = G' \oplus \mathcal{N}((G_0)_{\ell,\ell}) \oplus G'' \quad \text{with} \quad G' \cong ((\mathbb{G}_m[p^\infty])^f, \quad G'' \cong (\mathbb{Q}_p/\mathbb{Z}_p)^s.$$

We apply (2.2) in order to have the Dieudonné module of $(G_0)_{\ell,\ell}$ in normal form, obtaining a basis for the Dieudonné module of $G_0$. Now $\Diamond := \Diamond(\rho)$; this is the parallelogram of integral points $(x, y)$ bounded by $0 \leq y < c$, and $y > x \geq y - d$. The universal deformation $\mathcal{D}$ of $\mathbb{D}(G_0)$ is given by a display-matrix

$$\begin{pmatrix} A + TC & B + TD \\ C & D \end{pmatrix},$$

where $T = (T_{r,s} \mid 1 \leq r \leq d < s \leq h)$. Here $T_{r,s} = (t_{r,s}, 0 \cdots)$ and $\mathcal{D} = \mathrm{Spec}(k[[t_{r,s}]])$. We write as in (2.8) these variables in a diagram, by putting

$$T_{r,s} \quad \text{on the spot} \quad (s - r, s - 1 - d) \in \Diamond = \Diamond(\rho);$$

i.e.,

$$\begin{array}{cccccccccc}
 & & & & & & 0 & \cdots & 0 & -1 \\
 & & & & & T_{d,h} & \cdot & \cdots & T_{1,h} & \\
 & & & . & \vdots & \vdots & \vdots & . & & \\
 & T_{d,d+2} & \cdots & T_{i,d+2} & \cdots & T_{2,d+2} & T_{1,d+2} & & & \\
T_{d,d+1} & \cdots & T_{i,d+1} & \cdots & \cdots & T_{1,d+1} & & & &
\end{array}$$

We see that

$$\mathcal{D} = \mathrm{Spec}(R) = \mathrm{Spec}(k[[z_{(x,y)} \mid (x, y) \in \Diamond]]), \qquad T_{r,s} = Z_{(s-r,s-1-d)}.$$

For any $\beta \succ \mathcal{N}(G_0)$,

$$R_\beta = \frac{k[[z_{(x,y)} \mid (x, y) \in \Diamond]]}{(z_{(x,y)} \quad \forall (x, y) \notin \Diamond(\beta))}.$$

*Claim.*

$$(\mathrm{Spec}(R_\beta) \subset \mathrm{Spec}(R)) \quad = \quad (V_\beta \subset \mathcal{D})\,.$$

Clearly, the claim proves all statements of the theorem.

Since $G_0$ is a direct sum of a local-étale, a local-local and an étale-local $p$-divisible group as above, we obtain $\mathcal{N}(G_0) = \gamma = f{\cdot}(0, 1) + \gamma' + s{\cdot}(1, 0)$, where $\gamma' = \mathcal{N}((G_0)_{\ell,\ell})$. Note that $\beta' \mapsto f{\cdot}(0, 1) + \beta' + s{\cdot}(1, 0)$ gives a bijection between on the one hand all $\beta' \succ \gamma'$ and on the other hand all $\beta \succ \gamma$. Note moreover, see [1, Th. 4.4], that deformations of $G_0$ are smoothly fibered over deformations of $(G_0)_{\ell,\ell}$, with precise information on parameters (the matrix of the display is in blocks). This shows that the theorem follows for $G_0$ if it is proved for $(G_0)_{\ell,\ell}$. Hence we are reduced to proving the theorem in case $G_0$ is supposed to be of local-local type.

We use (2.8) in case of the normal form of the matrix of $\mathbb{D}(G_0)$ over $k$; we concluded that the entries $a_{i,j}$, with $1 \leq i \leq d$, and $d + 1 \leq j \leq h$ which are nonzero are all situated in $\Diamond(\mathcal{N}(G_0))$. For any integral domain $B$ which is a quotient of $R = k[[z_{(x,y)} \mid (x,y) \in \Diamond]]$ we can apply (2.6) and (2.8) to its field of fractions. This shows that $R \to B$ factors through $R \to R_\beta$ if and only if the generic fiber over $B$ has Newton polygon equal to or above $\beta$. This proves our claim in case $G_0$ is local-local. Hence this finishes the proof of the theorem. □

(3.3) We fix an integer $g$. For every *symmetric* Newton polygon $\xi$ of height $2g$,

$$\triangle(\xi) = \{(x,y) \in \mathbb{Z} \times \mathbb{Z} \mid y < g, \quad y < x \leq g, \quad (x,y) \prec \xi\},$$

and

$$\mathrm{sdim}(\xi) = \#(\triangle(\xi)).$$

Consider a $p$-divisible group $G_0$ over $k$ of dimension $g$ with a principal quasi-polarization. Now, $\mathcal{N}(G_0) = \gamma$; this is a symmetric Newton polygon. Now, $\mathcal{D} = \mathrm{Def}(G_0, \lambda)$ for the universal deformation space. For every symmetric Newton polygon $\xi$ with $\xi \succ \gamma$ we define $W_\xi \subset \mathcal{D}$ as the maximal closed, reduced formal subscheme carrying all fibers with Newton polygon equal to or above $\xi$; this space exists by Grothendieck-Katz; see [9, Th. 2.3.1, p. 143]. Note that $W_\rho = \mathcal{D}$, where $\rho = g\cdot((1,0) + (0,1))$.

(3.4) Theorem (NP-strata for principally quasi-polarized formal groups). *Suppose* $a(G_0) \leq 1$. *For every symmetric* $\xi \succ \gamma := \mathcal{N}(G_0)$, $\dim(W_\xi) = \mathrm{sdim}(\xi)$. *The strata* $W_\xi$ *are nested as given by the partial ordering on symmetric Newton polygons*; *i.e.*,

$$W_\xi \subset W_\delta \Longleftrightarrow \triangle(\xi) \subset \triangle(\delta) \Longleftrightarrow \xi \prec \delta.$$

*Generically on* $W_\xi$ *the fibers have Newton polygon equal to* $\xi$. *There is a coordinate system on* $\mathcal{D}$ *in which all* $W_\xi$ *are given by linear equations.*

(3.5) Corollary. *There exists a principally polarized abelian variety* $(X_0, \lambda_0)$ *over* $k$. *Strata in* $\mathrm{Def}(X_0, \lambda_0)$ *according to Newton polygons are exactly as in* (3.4). *In particular*, *the fiber above the generic point of* $W_\xi$ *is a principally polarized abelian scheme over* $\mathrm{Spec}(B_\xi)$ *having Newton polygon equal to* $\xi$ (*for* $B_\xi$, *see the proof of* (3.4) *below*; *for the notion* "*generic point of* $W_\xi$" *see the proof of* (3.5) *below*).

*Proof.* We write $(X_0, \lambda_0)[p^\infty] =: (G_0, \lambda_0)$. By Serre-Tate theory, see [8, §1], the formal deformation spaces of $(X_0, \lambda_0)$ and of $(G_0, \lambda_0)$ are canonically isomorphic, say $(\mathcal{X}, \lambda) \to \mathrm{Spf}(R)$ and $(\mathcal{G}, \lambda) \to \mathrm{Spf}(R)$ and $(\mathcal{X}, \lambda)[p^\infty] \cong (\mathcal{G}, \lambda)$. By Chow-Grothendieck, see [4, $\mathrm{III}^1$.5.4] (this is also called a theorem of "GAGA-type"), the formal polarized abelian scheme is algebraizable, and we obtain the universal deformation as a polarized abelian scheme over Spec$(R)$. Next, we can consider the generic point of $W_\xi \subset \mathrm{Spec}(R)$. Hence the Newton polygon of fibers can be read off from the fibers in $(\mathcal{G}, \lambda) \to \mathrm{Spec}(R)$. This proves that (3.5) follows from (3.4). $\square$

*Proof of* (3.4). The proof of this theorem is analogous to the proof of (3.2). We use the diagram

$$\begin{array}{cccccccc} & & & & & & -1 & \\ & & & X_{g,g} & \cdots & X_{1,g} & & \\ & & . & \vdots & & & & \\ 1 & X_{g,1} & \cdots & X_{1,1} & & & & . \end{array}$$

Here $X_{i,j},\ \ 1 \leq i, j \leq g$, is written on the place with coordinates $(g-i+j, j-1)$. We use the ring

$$B := \frac{k[[x_{i,j}; 1 \leq i, j \leq g]]}{(x_{k\ell} - x_{\ell k})}, \quad x_{i,j} = z_{(g-i+j,j-1)}, \quad (g-i+j, j-1) \in \triangle.$$

Note that $B = k[[x_{i,j} \mid 1 \leq i \leq j \leq g]] = k[[z_{x,y} \mid (x,y) \in \triangle]]$. For a symmetric $\xi$ with $\xi \succ \mathcal{N}(X_0)$ we consider

$$B_\xi = \frac{k[[t_{i,j}; 1 \leq i, j \leq g]]}{(t_{k\ell} - t_{\ell k}, \quad \text{and} \quad z_{(x,y)} \ \ \forall (x,y) \notin \triangle(\xi))}.$$

With this notation, applying (2.6) and (2.8), we finish the proof of (3.4) as we did in the proof of (3.2) above. $\square$

(3.6) *A remark on numbering.* In the unpolarized case, see (3.2), where there is a deformation matrix $(T_{r,s} \mid 1 \leq r \leq d < s \leq h)$, we put $T_{r,s}$ in the "NP diagram space" on the spot $(s-r, s-1-d)$. In the polarized case, see (3.4), we have a square matrix; according to notation used in case of the Riemann symmetry condition we write $(X_{i,j} \mid 1 \leq i, j \leq d = g)$, and put $X_{i,j}$ on the spot $(g-i+j, j-1)$, with $d = g = c$. Up to this change in numbering in the indices these amount to the same when methods concern the same cases: $T_{r,s} = X_{r,s-d}$. Note that we work with entries in the matix $(a)$ just below the diagonal, and obtain deformations starting from the variables $t_{r,s}$ producing elements "$a_{r,s-1} = t_{r,s}$", which cause the shifts in indices between (2.8) and (3.2).

## 4. Where to start

(4.1) LEMMA. *Given $d \in \mathbb{Z}_{>0}$, $c \in \mathbb{Z}_{>0}$ and a prime number $p$, there exists a field $K$ of characteristic $p$, and a formal group $G$ over $K$ of height $h = d + c$, of dimension $d$ with $a(G) = 1$ such that its Newton polygon $\mathcal{N}(G)$ is the straight line connecting $(0,0)$ with $(h,c)$, i.e. it has $h$ slopes each equal to $c/h$.*

(4.2) Suppose we are given $h, d \in \mathbb{Z}_{>0}$ with $d < h$. Consider the matrix (of size $h \times h$ with left hand upper corner a block of size $d \times d$):

$$\begin{pmatrix}
0 & 0 & \cdots & 0 & 0 & 0 & \cdot & \cdot & 0 & -p \\
1 & 0 & \cdots & 0 & 0 & 0 & & \cdots & \cdot & 0 \\
0 & 1 & \ddots & 0 & 0 & \cdot & & & \cdot & \cdot \\
\vdots & \vdots & \ddots & \ddots & \vdots & \vdots & & 0 & \vdots & \vdots \\
0 & 0 & \cdots & 1 & 0 & 0 & \cdot & \cdots & \cdot & 0 \\
0 & \cdots & \cdots & 0 & 1 & 0 & \cdot & \cdots & 0 & 0 \\
0 & \cdots & & \cdots & 0 & p & 0 & \cdots & \cdot & 0 \\
0 & \cdots & & \cdots & 0 & 0 & p & \ddots & \cdot & 0 \\
\vdots & & & & \vdots & \vdots & & \ddots & 0 & 0 \\
0 & \cdots & & \cdots & 0 & 0 & 0 & \cdots & p & 0
\end{pmatrix}.$$

We say this matrix is in *cyclic normal form* of height $h$ and dimension $d$.

*Proof of* (4.1). Consider (say over $K = \mathbb{F}_p$) the display given by the matrix above. Clearly this defines a Dieudonné module $M$ and hence (as long as $K$ is perfect) a $p$-divisible formal group $G$ with $\mathbb{D}(G) = M$. This matrix gives $F$ on a basis for $M$. We see that $F^h{\cdot}X_1 = -p^c{\cdot}X_1$. Hence $\mathcal{N}(G)$ is a straight line, e.g. apply (2.7). From the matrix we see that the Hasse-Witt matrix of $G$ has rank equal to $d-1$; hence $a(G) = 1$. $\square$

(4.3) LEMMA. *Given $g \in \mathbb{Z}_{>0}$ and a prime number $p$ there exists*:

(1) *A field $K$ of characteristic $p$, and a principally quasi-polarized supersingular formal group $(G, \lambda)$ over $K$ of dimension equal to $g$ with $a(G) = 1$.*

(2) *A field $k$ of characteristic $p$, and a principally polarized supersingular abelian variety $(X, \lambda)$ over $k$ of dimension equal to $g$ with $a(X) = 1$.*

*Proof.* Consider a matrix in cyclic normal form with $h = 2g$, and $d = g$. This defines a supersingular formal group $G$ of dimension $g$, height $h$, with $a(G) = 1$. We give a nondegenerate skew form on $M = \mathbb{D}(G)$ by requiring the basis to be symplectic; we are going to show this is possible. We compute the action of $V$:

$$V(X_{i+1}) = pX_i, \quad VX_1 = -Y_g, \quad VY_1 = pX_g, \quad VY_{j+1} = Y_j.$$

Now, $F$ and $V$ must respect the pairing on this symplectic base; i.e., we have to prove that for $a, b \in M$,

$$\langle Fa, b\rangle \quad = \quad \langle a, Vb\rangle^{\sigma}.$$

It suffices to show this directly on base vectors; let us show the essential steps (the others being obvious):

$$-1 = \langle Y_1, X_1\rangle = \langle FX_d, X_1\rangle \stackrel{?}{=} \langle X_g, VX_1 = -Y_g\rangle^{\sigma} = -1,$$

and analogously $\langle FY_d, Y_1\rangle = -p = \langle Y_d, VY_1\rangle^{\sigma}$. This proves the first statement.

In order to prove (2), we have to see that $(G, \lambda)$ is algebraizable, i.e. comes from a principally polarized abelian variety. Using [21, §2] there is exists (canonically) a polarized flag-type quotient $H \to H/N = G$, $(H, \mu') \to (G, \lambda)$ with $a(H) = g$. We know that $a(H) = g$ implies that $H \otimes k \cong E[p^{\infty}]^g$, where $E$ is a supersingular elliptic curve (isomorphism, say, over an algebraic closure $k$ of $K$). We see that $\mu'$ on $H$ comes from a polarization $\mu$ on $E^g$ with $\deg(\mu') = \deg(\mu)$ (e.g. use [13, Prop. 6.1]). Hence $(H, \mu')$ is algebraizable: $(E^g, \mu)[p^{\infty}] \cong (H, \mu') \otimes k$ and thus the quotient $(G, \lambda) \otimes k$ is algebraizable. This ends a construction which shows (2). □

(4.4) *Remark.* For a proof of (4.3) we could also refer to [13, (4.9)] (where it is proved that every component of the supersingular locus in the principally polarized case has generically $a = 1$). However that is a much more involved result than just the lemma above.

## 5. A conjecture by Manin

(5.1) A CONJECTURE BY MANIN (see [14, p. 76]). *For any prime number $p$ and any symmetric Newton polygon $\xi$ there exists an abelian variety $X$ over a field in characteristic $p$ such that $\mathcal{N}(X) = \xi$.*

This is the converse of (1.7).

(5.2) This conjecture was proved in the Honda-Serre-Tate theory; see [29, p. 98].

Using that theory we can actually prove somewhat more; we know that a supersingular abelian variety (of dimension at least 2) is not absolutely simple; however this is about the only general exception: for a symmetric Newton polygon which is not supersingular, there exists an absolutely simple abelian variety in characteristic $p$ having this Newton polygon; see [12]. Can we prove the result of that paper using the deformation theory as discussed here?

Once the conjecture by Manin is proved, we conclude that actually there exists an abelian variety defined over a finite field with the desired Newton polygon.

(5.3) *A proof of the Manin conjecture.* Suppose there is a symmetric Newton polygon $\xi$ of height $2g$. Choose $(X_0, \lambda_0)$ as in (4.3)(2) above; in particular $a(X_0) \leq 1$ and $\lambda_0$ is a principal polarization. Apply (3.5). We conlude that the generic fiber over $W_\xi \subset \mathrm{Def}(X_0, \lambda_0)$ is a (principally) polarized abelian variety with Newton polygon equal to $\xi$. □

(5.4) Remark. Actually we proved that for any $\xi$ there exists a principally polarized $(X, \lambda)$ with $\mathcal{N}(X) = \xi$ and $a(X) \leq 1$.

(5.5) Let us analyze the essential step, where we made a formal group algebraic. Starting with supersingular formal groups, we know that the $a$-number is maximal if and only if the formal group is isomorphic with $(G_{1,1})^g$ (see [23]); hence this can be algebraized (note that $E[p^\infty] \cong G_{1,1} \otimes \overline{\mathbb{F}_p}$ for a supersingular elliptic curve $E$). We can algebraize a principally quasi-polarized supersingular formal group with $a(G_0) = 1$ (use polarized flag-type quotients). Then we apply the deformation theory (say of quasi-polarized formal groups), obtaining the Newton polygon desired, and then the theory of Serre-Tate, see [8, §1]. Then the Chow-Grothendieck algebraization, see [4, $\mathrm{III}^1$.5.4], allows us to algebraize the family; it produces an abelian variety as required in the Manin conjecture.

(5.6) By methods just explained a weak form of the conjecture by Grothendieck follows:

Weak form. *Let $\beta \succ \gamma$ be two Newton polygons belonging to a height $h$ and a dimension $d$ which are comparable. There exists a $p$-divisible group $\mathcal{G} \to S$ over an integral formal scheme $S$ in characteristic $p$, such that the generic fiber $\mathcal{G}_\eta$ has $\mathcal{N}(\mathcal{G}_\eta) = \beta$ and the closed special fiber $G_0$ has $\mathcal{N}(\mathcal{G}_0) = \gamma$.*

Weak form, AV. *Suppose there are two* symmetric *Newton polygons, $\xi \succ \gamma$. There exists a specialization of polarized abelian varieties $(X_\eta, \lambda_\eta) \subset (\mathcal{X}, \lambda) \supset (X_0, \lambda_0)$ having $\xi$ and $\gamma$ as Newton polygons for the generic an special fiber.*

This proves the following (see Koblitz, [10, p. 215]): *Suppose there exists a sequence $\xi_1 \prec \cdots \prec \xi_n$ of comparable symmetric Newton polygons. Then there exists a sequence of specializations having these Newton polygons in each step.*

## 6. A conjecture by Grothendieck

(6.1) In [3, Appendix], we find a letter of Grothendieck to Barsotti, and on page 150 we read: "$\cdots$ The wishful conjecture I have in mind now is the following: the necessary conditions $\cdots$ that $G'$ be a specialization of $G$ are also sufficient. In other words, starting with a BT group $G_0 = G'$, taking its formal modular deformation $\cdots$ we want to know if every sequence of rational numbers satisfying $\cdots$ these numbers occur as the sequence of slopes of a fiber of $G$ as some point of $S$."

We say that $G_0$ is a *specialization* of $\mathcal{G}_\eta$ if there exists an integral local scheme $S$ and a $p$-divisible group $\mathcal{G} \to S$ with $G_0$ as closed fibre and $\mathcal{G}_\eta$ as generic fibre; in this case we shall write $\mathcal{G}_\eta \rightsquigarrow G_0$. We use analogous notation for polarized abelian varieties and quasi-polarized $p$-divisible groups.

(6.2) Theorem (a conjecture by Grothendieck). *Given a $p$-divisible group $G_0$ and Newton polygons $\mathcal{N}(G_0) =: \gamma \prec \beta$, assume $a(G_0) \leq 1$. Then there exists a specialization $\mathcal{G}_\eta \rightsquigarrow G_0$ with $\beta = \mathcal{N}(\mathcal{G}_\eta)$.*

(6.3) Theorem (an analogue of the conjecture by Grothendieck). (a) *Suppose there are a principally quasi-polarized $p$-divisible group $(G_0, \lambda_0)$ and symmetric Newton polygons $\mathcal{N}(G_0) =: \gamma \prec \xi$. Assume $a(G_0) \leq 1$. There exists a specialization $(\mathcal{G}_\eta, \mu) \rightsquigarrow (G_0, \lambda_0)$ with $\xi = \mathcal{N}(\mathcal{G}_\eta)$.*

(b) *Given a principally polarized abelian variety $(X_0, \lambda_0)$ and symmetric Newton polygons $\mathcal{N}(X_0) =: \gamma \prec \xi$, assume $a(X_0) \leq 1$. Then there exists a specialization $(\mathcal{X}_\eta, \mu) \rightsquigarrow (X_0, \lambda_0)$ with $\xi = \mathcal{N}(\mathcal{X}_\eta)$.*

*Proof of* (6.2). Applying (2.2) to $\mathbb{D}(G_0)$, we use deformation theory and its methods, as developed in Section 4, to obtain a generic fiber with the desired Newton polygon; by (3.2) any $\beta \succ \mathcal{N}(G_0)$ is realized in $\mathcal{D} = \mathrm{Def}(G_0)$. This proves (6.2).

*Proof of* (6.3). We apply (2.3) to $\mathbb{D}(G_0)$, respectively to $\mathbb{D}(X_0[p^\infty])$. By (3.4) and (3.5), a principally polarized $(X_0, \lambda_0)$ can be deformed to a principally polarized abelian variety with a given symmetric $\xi \succ \mathcal{N}(X_0)$. This proves all existence results in (6.3). $\square$

(6.4) *Remark. The analogue of the conjecture by Grothendieck does not hold for arbitrary polarized abelian varieties.* This is shown by the following:

*Example* (see [7, Remark 6.10]). Fix a prime number $p$, consider abelian varieties of dimension 3 with a polarization of degree $p^6$. In that 6-dimensional moduli space $\mathcal{A}_{3,p^3} \otimes \mathbb{F}_p$ (of course) the locus where the $p$-rank is zero, has dimension 3; see [18, Th. 4.1]. It can be proved that the supersingular locus in this moduli space has a component of dimension equal to three; see [21, Cor. 3.4]. Hence we conclude there exists a polarized abelian variety $(X_0, \lambda_0)$ of dimension 3, supersingular, hence $\mathcal{N}(X_0) = \sigma = 3{\cdot}(1,1)$ with $\mathrm{degree}(\lambda_0) = p^3$, such that every deformation of this polarized abelian variety either is supersingular or has positive $p$-rank; thus no such deformation will produce the Newton polygon $\gamma = (2,1) + (1,2)$. This ends the description of the example, and the claim that (6.2) does not hold for arbitrary polarized abelian varieties is proved.

We expect an example $(G_0, \lambda_0)$, which cannot be deformed to some $\xi \succ \gamma = \mathcal{N}(G_0)$ as just explained, to be available for every symmetric Newton polygon $\gamma$ with $f(\gamma) \leq g-2$, i.e. which allows in its isogeny class a formal group with $a \geq 2$, and a carefully chosen inseparable polarization.

(6.5) *Remark.* In the previous example we know that $a(X_0) = 1$ (by [18, Th. 4.1]). Hence for any Newton polygon $\gamma$ with $h = 6$, $d = 3$ the $p$-divisible group $X_0[p^\infty]$ (no quasi-polarization considered) can be deformed to a $p$-divisible group with Newton polygon equal to $\gamma$. We see the curious fact that a deformation of a given Newton polygon as $p$-divisible groups does exist, but as (nonprincipally) *polarized* $p$-divisible groups does not exist in this case.

(6.6) *Remark.* The conjecture by Grothendieck in its general form, i.e. (6.2) without assuming anything about $a(G_0)$, and the analogue for the (quasi-) principally polarized case, i.e. (6.3) without assuming anything about $a(G_0)$, have been proved; see [6] and [25].

It follows that the dimension (of every component) of $W_\xi \subset \mathcal{A}_{g,1} \otimes \mathbb{F}_p$ equals $\mathrm{sdim}(\xi)$, as announced in [24]; this was conjectured for the supersingular Newton polygon $\sigma$ in [21], and proved for $W_\sigma = \mathcal{S}_{g,1}$ in [13]. Note that for a Newton polygon stratum for $\xi$ for *nonseparably polarized* abelian varieties the dimension count can give an answer different from $\mathrm{sdim}(\xi)$.

UTRECHT UNIVERSITY, UTRECHT, THE NETHERLANDS
*E-mail address*: oort@math.uu.nl